# THE OSCILLATORY DISTRIBUTION OF DISTANCES IN RANDOM TRIES

By Costas A. Christophi and Hosam M. Mahmoud

*George Washington University*

We investigate $\Delta_n$, the distance between randomly selected pairs of nodes among $n$ keys in a random *trie*, which is a kind of digital tree. Analytical techniques, such as the Mellin transform and an excursion between poissonization and depoissonization, capture small fluctuations in the mean and variance of these random distances. The mean increases logarithmically in the number of keys, but curiously enough the variance remains $O(1)$, as $n \to \infty$. It is demonstrated that the centered random variable $\Delta_n^* = \Delta_n - \lfloor 2 \log_2 n \rfloor$ does not have a limit distribution, but rather oscillates between two distributions.

**1. Introduction.** In this computer age that we are living in, digital data, which are represented and accessed via their composition into digits, are quite ubiquitous in science and technology. They are prevalent in computer files, telecommunication signals, DNA and so forth. Digital searching deals with the storage of information in a computer's memory and its fast recovery [9]. Various forms of digital trees are known to be super efficient for this purpose. The subject of this investigation is the *trie*, one popular such digital technology, initially proposed independently by De La Briandais [2] and Fredkin [4] for information re*trie*val.

In informatics, distances between nodes in a random combinatorial object are of prime interest, because they are indicative of the speed of communication within the structure. These distances have applications in many other fields too. For example, the collective sum of all such distances in the graph underlying a molecule is known in chemistry as the *Wiener index* (see [5] and [19]).

In this paper we look at the distances between distinct pairs of keys in random binary tries. The mean and variance and ultimately the asymptotic distribution (via its moment generating function) are derived by analytical methods involving the use of poissonization, as a mathematical transform, and of depoissonization, as an asymptotic inverse transform. Although the chief interest lies in studying the random tree for a fixed population of $n$ keys, the recurrence equations involved are rather unwieldy. If a Poisson number of keys is assumed instead, the functional equations involved can asymptotically be solved by the Mellin transform and its inverse. It can then be justified that the solution is a good approximation (with quantifiable small errors) for the fixed-population problem, when the Poisson parameter is taken to be $n$, as $n \to \infty$.

Analogous studies have been conducted on binary search trees (a different model of random trees), as, for example, the study of the Wiener index [14], and the study of random distances [13] and its generalization to spanning trees of subsets of nodes [15]. However, the methodology for binary search trees is essentially different from the analytic tool-kit we employ here for tries.

Assume a standard Bernoulli probability model on data. Let $\Delta_n$ be the distance between two randomly chosen keys in a trie. Two of our results concern the mean and variance, and the periodic fluctuations therein (the notation lg stands for $\log_2$):

$$\mathbf{E}[\Delta_n] = 2\lg n + \eta(\lg n) - \frac{\ln 2 - 2\gamma}{\ln 2} + O\Big(\frac{1}{n^{0.4999}}\Big),$$

$$\mathbf{Var}[\Delta_n] = \frac{2\pi^2 + 19\ln^2 2 - \alpha}{3\ln^2 2} + \xi(\lg n) + O\Big(\frac{1}{n^{0.4999}}\Big),$$

where $\alpha \approx 7.227113$ is a constant explicitly found as a numerical series, and $\eta(\cdot)$ and $\xi(\cdot)$ are bounded functions of miniscule magnitude which are periodic in their argument. The error terms include additional small periodic functions as well. Note that the variance is $O(1)$, a highly desirable feature in a data structure; this property indicates robustness of random digital trees in spite of the large variation in the data.

In view of the $O(1)$ variance, one expects that if a nontrivial limit distribution were to exist, it would require only centering, but no scaling. Our main result is that a centered integer version of $\Delta_n$ does not converge in distribution to any random variable. Rather its distribution function assumes the form of a discrete staircase distribution which oscillates between two distinct fixed discrete distribution functions.

The main results have been sketched. The rest of the paper is organized in sections as follows. In Section 2 the definition of a trie is given and illustrated by an example, and the associated trie terminology is made precise. The probability distribution assumed on the data is also discussed. At the end of Section 2, notation that is used as a working language throughout is explained. The moments are discussed in Section 3 and the moment generating function is derived in Section 3.1. The poissonization-Mellin-inverse Mellin-depoissonization program is illustrated in this context. The residue calculation is taken up in Section 3.2. Some of the very lengthy calculations are relegated to Appendixes A and B. The first two moments are pumped out from the Mellin transform of the moment generating function. This is done thoroughly for the mean in Section 3.3 and only sketched for the variance in Section 3.4. The oscillating nature of the moment generating function elicits the nonexistence of a limit distribution. This is shown in Section 3.5 by demonstrating explicit fluctuations in the moment generating function. Section 4 concludes the paper with an interpretation of the results and a discussion of their scope.

**2. Tries.** A *trie* is a digital tree structure consisting of *internal nodes* that have one or two children, and *leaf nodes* that hold data items (commonly called *keys*). We shall assume our data to be binary (although our methods should work well for larger alphabets too). Suppose we have $n \ge 0$ keys given in their dyadic representation. Data are scaled to be in the interval $[0,1]$. The trie grows by the progressive insertion of keys. If $n=0$, nothing needs to be done; the insertion algorithm terminates. If $n=1$, a leaf is allocated for the single key given. If $n \ge 2$, an internal node is allocated as a *root* of the tree; keys with 0 as their most significant bit go to the left subtree, and keys with 1 to the right. Subsequently, in the subtrees the insertion algorithm is applied recursively, using the $(\ell+1)$st most significant bit of the key for branching at level $\ell$. At termination, each key is stored in a leaf by itself, and the root-to-leaf paths in the tree correspond to minimal prefixes sufficient to distinguish the keys. As an illustration, suppose that $n=5$ and that the data are the following:

$$X_1 = 0.00111\dots,$$
$$X_2 = 0.11011\dots,$$
$$X_3 = 0.00011\dots,$$
$$X_4 = 0.01010\dots,$$
$$X_5 = 0.11111\dots.$$

The algorithm guides $X_1, X_3$ and $X_4$ to the left subtree, and it guides $X_2$ and $X_5$ to the right one. Then, in the left subtree $X_1$ and $X_3$ go to the left but $X_4$ goes to the right, whereas in the right subtree, both $X_2$ and $X_5$ go to the right. The subtree containing $X_1$ and $X_3$ is developed further into a left subtree with $X_3$ being a leaf node as it is the only node in the subtree, and a right subtree with $X_1$ contained in a leaf node as it is the only node in the subtree, and so on. The resulting trie is depicted in Figure 1, with bullets indicating internal nodes and squares indicating leaf nodes, where the keys are stored.

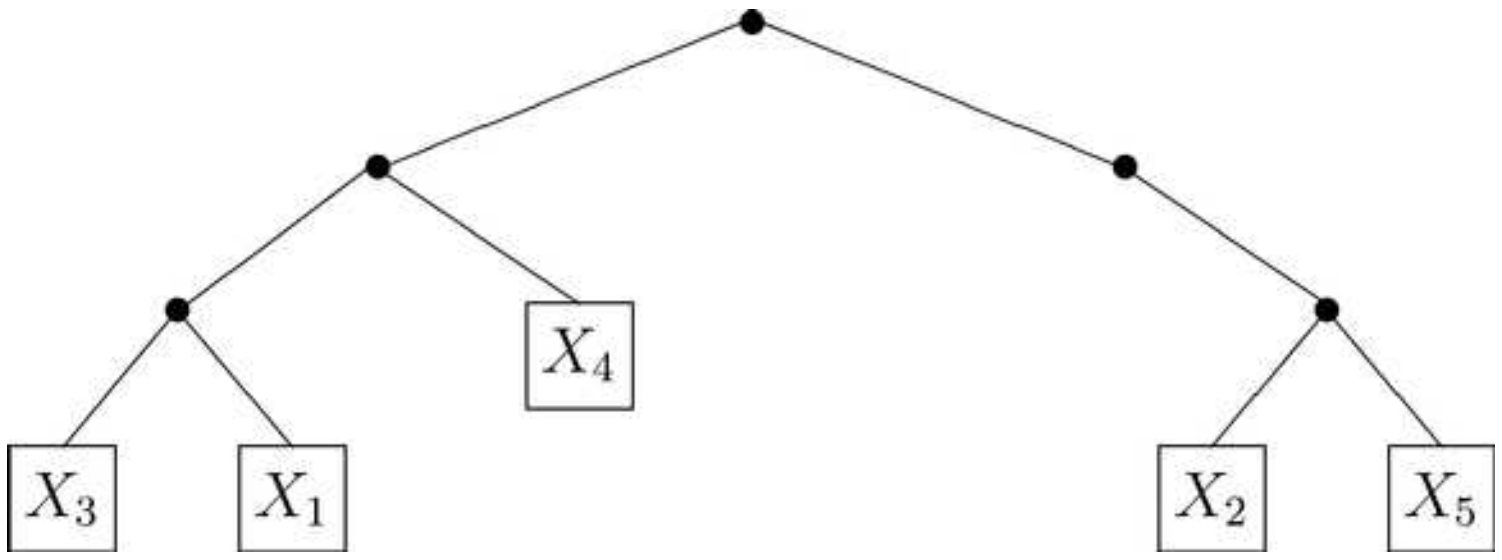

FIG. 1. *A trie with 5 keys.*

We assume a standard model of randomness, according to which keys are independent and of infinite length; within each the bits are independent and equiprobable. Thus, each key can be viewed as an infinite sequence of Bernoulli trials. This probability model is often called the *unbiased Bernoulli model*.

We define $\Delta_n$ as the distance (i.e., the number of tree edges) between two randomly selected keys in a random trie of size $n$, with *random* meaning that all $\binom{n}{2}$ pairs of keys are equally likely choices. For example, in Figure 1 the distance between $X_2$ and $X_4$ is 5. A recursive formulation for $\Delta_n$ goes via $\delta_n$, the depth of a randomly selected key in a random trie of size $n$.

Furthermore, the Mellin transform of a function $f(x)$ is

$$\int_0^\infty f(x)x^{s-1}\,dx,$$

and will be denoted by $\mathbf{Mellin}\{f(x);s\}$ or, interchangeably, by $f^*(s)$. The Mellin transform usually exists in vertical strips in the $s$ complex plane of the form

$$a < \Re s < b,$$

for real numbers $a < b$. We shall denote such a domain of existence by $\langle a, b\rangle$. The function $f(x)$ can be recovered from its transform by a line integral

$$f(x) = \frac{1}{2\pi i}\int_{c-i\infty}^{c+i\infty} f^*(s)x^{-s}\,ds,$$

for any $c \in (a, b)$.

**3. Moments of the random distance.** Let $L_n$ and $R_n$ be, respectively, the number of keys residing in the left and right subtrees, among the $n$ keys received by the tree (so, $L_n + R_n = n$). In view of the unbiased Bernoulli model, $L_n$ and $R_n$ are both distributed like $\mathrm{Bin}(n, \frac{1}{2})$, a binomially distributed random variable on $n$ independent trials and rate of success $\frac{1}{2}$ per trial. Because of the probability distribution of the keys and the recursive action of the insertion algorithm, each of the left and the right subtree preserves the probabilistic structure of the trie. The only difference is that they are *random* tries constructed on $L_n$ and $R_n$ keys, respectively, instead of $n$.

Given $L_n$, $\Delta_n$ can be $\Delta_{L_n}$ with probability $\binom{L_n}{2}/\binom{n}{2}$ when both keys come from the left subtree, $\tilde{\Delta}_{R_n}$ with probability $\binom{R_n}{2}/\binom{n}{2}$ when both keys come from the right subtree or $(\delta_{L_n}+1)+(\tilde{\delta}_{R_n}+1)$ with probability $L_nR_n/\binom{n}{2}$ when the keys come from different subtrees [here a tilded random variable indicates an independent copy (with the same distribution) of that random variable]. Hence, we have the conditional distribution

$$\Delta_n | L_n = \begin{cases} \Delta_{L_n}, & \text{with probability } \dfrac{\binom{L_n}{2}}{\binom{n}{2}}, \\ \tilde{\Delta}_{R_n}, & \text{with probability } \dfrac{\binom{R_n}{2}}{\binom{n}{2}}, \\ (\delta_{L_n}+1)+(\tilde{\delta}_{R_n}+1), & \text{with probability } \dfrac{L_nR_n}{\binom{n}{2}}, \end{cases} \tag{1}$$

where the vector $(\tilde{\Delta}_1, \ldots, \tilde{\Delta}_n)$ denotes a copy of random distances independent of the vector $(\Delta_1, \ldots, \Delta_n)$; likewise, $(\tilde{\delta}_1, \ldots, \tilde{\delta}_n)$ and $(\delta_1, \ldots, \delta_n)$ are vectors of independent copies of random depths. It is important to note that $\Delta_{L_n}$ and $\tilde{\Delta}_{R_n}$ are *dependent* through the dependency of $L_n$ and $R_n$, but given the value of $L_n$ (and, hence, $R_n$), the two are *conditionally independent*; the same applies to $\delta_{L_n}$ and $\tilde{\delta}_{R_n}$.

3.1. *The moment generating function.* From the conditional recursion (1), we have a recurrence for the conditional moment generating function:

$$\binom{n}{2}\mathbf{E}[e^{\Delta_n t}|L_n] = e^{\Delta_{L_n} t}\binom{L_n}{2} + e^{\tilde{\Delta}_{R_n} t}\binom{R_n}{2} + e^{((\delta_{L_n}+1)+(\tilde{\delta}_{R_n}+1))t}L_nR_n.$$

Taking expectation of both sides, we essentially get a recurrence for the moment generating function:

$$\binom{n}{2}\mathbf{E}[e^{\Delta_n t}] = \mathbf{E}\left[e^{\Delta_{L_n} t}\binom{L_n}{2}\right] + \mathbf{E}\left[e^{\tilde{\Delta}_{R_n} t}\binom{R_n}{2}\right] + \mathbf{E}[e^{((\delta_{L_n}+1)+(\tilde{\delta}_{R_n}+1))t}L_nR_n].$$

We observe the symmetry of the left and the right subtrees, and the fact that $L_n$ is distributed like $\mathrm{Bin}(n, \frac{1}{2})$. After multiplying by $z^n$ and summing over all possible values of $n$, we get

$$\sum_{n=0}^{\infty} z^n \frac{\binom{n}{2}\phi_{\Delta_n}(t)}{n!} = 2\sum_{n=0}^{\infty}\sum_{\ell=0}^{n} z^n \frac{\binom{\ell}{2}\phi_{\Delta_\ell}(t)}{\ell!\,(n-\ell)!2^n} + \sum_{n=0}^{\infty}\sum_{\ell=0}^{n} z^n \mathbf{E}\Big[e^{(\delta_\ell+1)t}\frac{\ell}{\ell!}\Big]\mathbf{E}\Big[e^{(\delta_{n-\ell}+1)t}\frac{n-\ell}{(n-\ell)!}\Big]\frac{1}{2^n},$$

where for any generic random variable $X$, $\phi_X(t) := \mathbf{E}[e^{tX}]$. This gives

$$\Phi(t,z) = 2\Phi\Big(t,\frac{z}{2}\Big)e^{z/2} + \Big\{\sum_{j=0}^{\infty}\mathbf{E}\Big[e^{(\delta_j+1)t}\frac{j}{j!}\Big]\Big(\frac{z}{2}\Big)^j\Big\}^2, \qquad (2)$$

where $\Phi(t,z) := \sum_{j=0}^{\infty} z^j \frac{\binom{j}{2}\phi_{\Delta_j}(t)}{j!}$.

Note that $e^{-z}\Phi(t,z)$ has a poissonization interpretation. For, if $N(z)$ is distributed like a Poisson random variable with mean $z$, then

$$\begin{aligned} e^{-z}\Phi(t,z) &= \sum_{j=0}^{\infty}\binom{j}{2}\phi_{\Delta_j}(t)P(N(z)=j) \\ &= \mathbf{E}\Big[\binom{N(z)}{2}e^{\Delta_{N(z)}t}\Big]. \end{aligned}$$

Averaging $\binom{N(z)}{2}e^{\Delta_{N(z)}t}$ is just like averaging $\binom{n}{2}e^{\Delta_n t}$, except that we assume a Poisson number of keys instead of a fixed population $n$. As we shall see, the advantage of poissonization is that the Poisson model is amenable to the formulation of functional equations that can be asymptotically solved by the Mellin transform techniques. The various fixed-population problems are averaged over cases with Poisson probabilities as weights. The fixed-population asymptotics can be well approximated by their depoissonized counterparts.

When we multiply (2) throughout by $e^{-z}$, we get the functional equation

$$e^{-z}\Phi(t,z) = 2e^{-z/2}\Phi\Big(t,\frac{z}{2}\Big) + F\Big(t,\frac{z}{2}\Big),$$

where $F(t,z) := [e^t e^{-z}\sum_{j=0}^{\infty}\frac{j\phi_{\delta_j}(t)}{j!}z^j]^2$.

We want to study convergence issues of moment generating functions in a neighborhood around $t=0$ (via the Mellin transform). The Mellin transform of the latter functional equation does not exist in such a neighborhood to recover the moment generating function. For example, at $t=0$, the latter equation becomes

$$\frac{z^2}{2} = \frac{z^2}{4} + \frac{z^2}{4},$$

and we cannot take the Mellin transform of the equation because $z^2$ does not have such a transform. Nevertheless, the same equation shifted down by $z^2/2$ has a Mellin transform. It turns out that for general $t$, the factor $e^{2t}\frac{z^2}{2}$ is a suitable shift. Introduce this shift to get

$$e^{-z}\Phi(t,z) - e^{2t}\frac{z^2}{2} = 2\left[e^{-z/2}\Phi\left(t,\frac{z}{2}\right) - e^{2t}\frac{(z/2)^2}{2}\right] - e^{2t}\frac{z^2}{4} + F\left(t,\frac{z}{2}\right),$$

that is,

$$P(t,z) = 2P\left(t,\frac{z}{2}\right) + F\left(t,\frac{z}{2}\right) - e^{2t}\frac{z^2}{4}, \qquad (3)$$

where $P(t,z) := e^{-z}\Phi(t,z) - e^{2t}\frac{z^2}{2}$.

LEMMA 1. *The Mellin transform of $P(t,z)$ (with respect to the variable $z$) is*

$$\begin{aligned} P^*(t,s) &= \frac{e^{2t}2^s(1-e^t)^2\Gamma(s+2)}{(1-2^{s+1})(1-e^{2t}2^{s+2})} \\ &\quad\times\left(-\frac{2^{s+3}e^t}{1-e^t2^{s+2}} - 2e^t\sum_{k=0}^{\infty} e^{tk}\left[1 - \frac{1}{(1+1/(2^{k+1}))^{s+2}}\right] + \frac{(1-2^{s+3})}{2^{s+2}}\right) \\ &\quad+ \frac{2^{s+1}e^{2t}(1-e^t)\Gamma(s+2)}{(1-2^{s+1})(1-e^t2^{s+2})}, \end{aligned}$$

*well defined for $\langle -3, -2 - \frac{2|t|}{\ln 2}\rangle$, for $|t| < \frac{1}{2}\ln 2$.*

PROOF. The bulk of this computation is in the handling of the function $F(t,z)$, which is a squared generating function. We relegate this lengthy computation to Appendixes A and B. □

Note that, for negative $t$, the strip $\langle -3, -2 - \frac{2t}{\ln 2}\rangle$ includes the entire strip $\langle -3, -2\rangle$, and, for all $|t| < \frac{1}{2}\ln 2$, the Mellin transform exists in $\langle -3, -2 - \frac{2|t|}{\ln 2}\rangle$. The inverse Mellin transform,

$$P(t,z) = \frac{1}{2\pi i}\int_{c-i\infty}^{c+i\infty} P^*(t,s)z^{-s}\,ds,$$

recovers the (shifted) poissonized moment generating function, if we take $c \in (-3, -2 - \frac{2|t|}{\ln 2})$, with $|t| < \frac{1}{2}\ln 2$. The transform $P^*(t,s)$ has simple poles at $s = -2$, at $s = -1 + \frac{2\pi ri}{\ln 2}$, for $r = 0, \pm 1, \ldots$, at $s = -2 + \frac{-2t+2\pi ji}{\ln 2}$, for $t \neq 0$ and $j = 0, \pm 1, \ldots$, at $s = -2 + \frac{-t+2\pi hi}{\ln 2}$, for $t \neq 0$ and $h = 0, \pm 1, \ldots$ and at $s = -3, -4, \ldots$.

We employ the method of "closing the box." (This method is discussed in [12] and [18].) In this method one takes the complex integration over the line $c-iM$ and $c+iM$, and then closes the box connecting the four corners $c \pm iM$ and $d \pm iM$, for an arbitrary $d > 0$. The number $M$ is chosen in such a way that no pole is crossed. For example, we can take it to be $\frac{\pi i(2m+1)}{\ln 2}$, for integer $m$. The residue theorem points out that

$$\lim_{M\to\infty} \oint P^*(t,s) z^{-s}\, ds = 2\pi i \sum \text{residues in } \langle c, d\rangle. \tag{4}$$

The contour integral can be written as

$$\oint P^*(t,s) z^{-s}\, ds = \int_{c-iM}^{d-iM} + \int_{d-iM}^{d+iM} + \int_{d+iM}^{c+iM} + \int_{c+iM}^{c-iM}.$$

As we let $m \to \infty$ (hence, $M \to \infty$), the line integrals at the top and bottom sides of the box approach 0, as the magnitude of the $\Gamma$ function decreases exponentially fast with its imaginary part. Moreover, the integral at the right side of the box introduces an error term of the order of $O(z^{-d})$. Hence, (4) gives

$$P(t,z) = O(z^{-d}) - \sum \text{residues in } \langle c, d\rangle.$$

The problem has now been reduced to residue computation.

3.2. *Residue computation.* We have the following residue calculations:

$$\underset{s=-1}{Res}[P^*(t,s)z^{-s}] = \frac{z(1-e^t)e^{2t}}{4(1-2e^{2t})\ln 2}$$

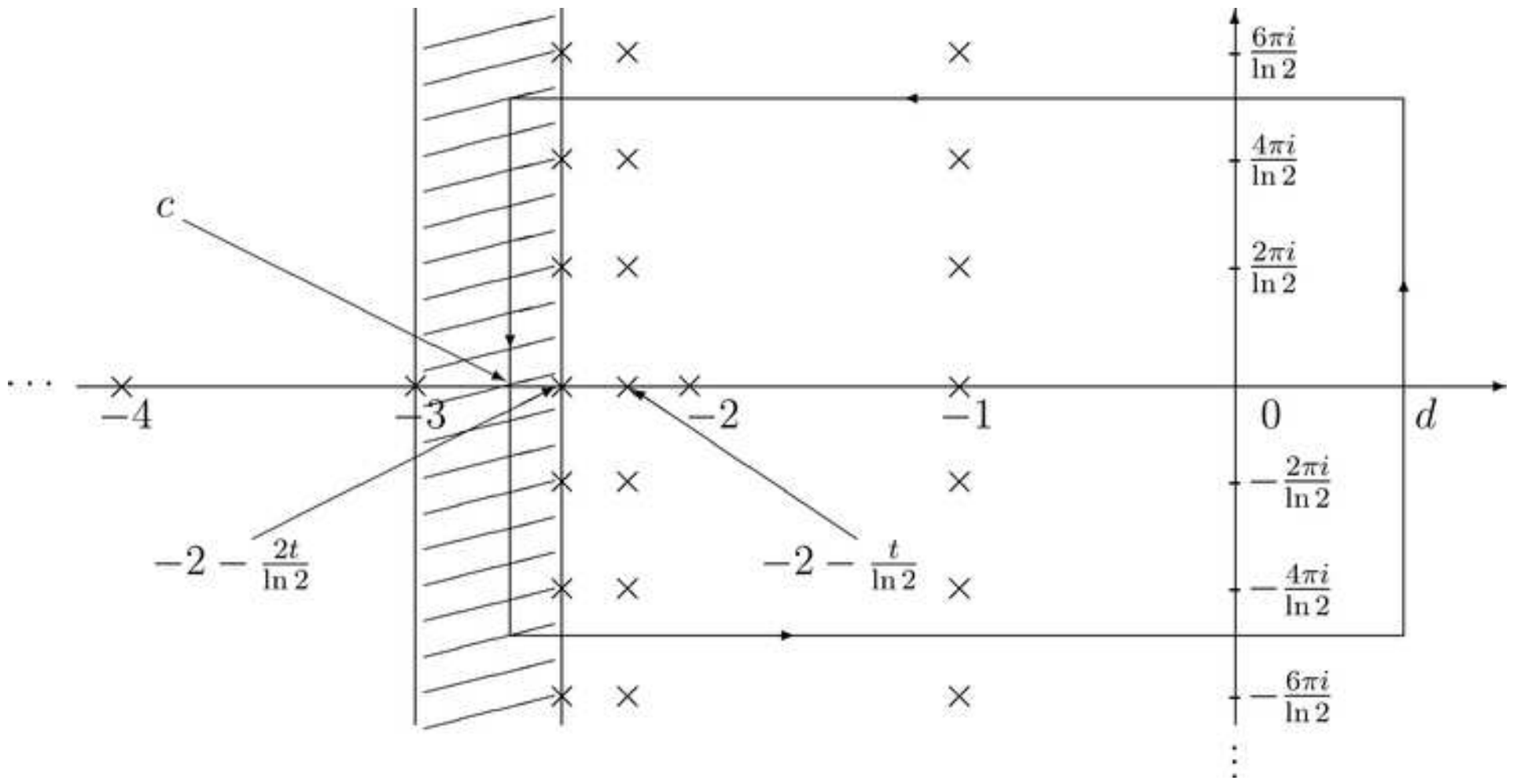

Fig. 2. *The poles of the Mellin transform $P^*(t,s)$; $\times$ indicates a simple pole. The shaded area is a domain free of poles, where $P^*(t,s)$ is well defined.*

$$\times \left[ 4e^t \sum_{k=0}^{\infty} \frac{e^{tk}(1-e^t)}{2^{k+1}+1} - 3e^t - 1 \right],$$

$$\operatorname*{Res}_{s=-2}[P^*(t,s)z^{-s}] = \tfrac{1}{2} z^2 e^{2t},$$

$$\operatorname*{Res}_{s=-2-t/\ln 2}[P^*(t,s)z^{-s}] = 0,$$

$$\operatorname*{Res}_{s=-2-2t/\ln 2}[P^*(t,s)z^{-s}] = \frac{z^{2+2t/\ln 2}\Gamma(-2t/\ln 2)e^{3t}(e^t-1)}{(1-2e^{2t})\ln 2} \times \Bigg( \frac{1}{2}(e^{2t}-e^t-2) - (e^t-1)\sum_{k=0}^{\infty} e^{tk}\left[1 - \left(1+\frac{1}{2^{k+1}}\right)^{2t/\ln 2}\right]\Bigg),$$

$$\operatorname*{Res}_{s=-1+2\pi ri/\ln 2}[P^*(t,s)z^{-s}] = (e^t-1)\frac{2\pi irz\Gamma(2\pi ir/\ln 2)e^{-2\pi ir \lg z}}{(1-2e^{2t})\ln^2 2} \times \Bigg( \frac{e^{2t}(e^t-2)}{1-2e^t} - e^{3t}(e^t-1) \times \sum_{k=0}^{\infty} e^{tk}\left[1 - \left(1+\frac{1}{2^{k+1}}\right)^{-1-2\pi ir/\ln 2}\right]\Bigg),$$

$$r \neq 0,$$

$$\operatorname*{Res}_{s=-2+(-t+2\pi hi)/\ln 2}[P^*(t,s)z^{-s}] = 0, \qquad h \neq 0,$$

and, for $j \neq 0$,

$$\operatorname*{Res}_{s=-2+(-2t+2\pi ji)/\ln 2}[P^*(t,s)z^{-s}] = \frac{z^{2+2t/\ln 2}e^{-2\pi ij \lg z}\Gamma((-2t+2\pi ij)/\ln 2)(e^t-1)}{(1-2e^{2t})\ln 2} \times \Bigg( \frac{e^{2t}-e^t-2}{2} - e^{3t}(e^t-1)\sum_{k=0}^{\infty} e^{tk}\left[1 - \left(1+\frac{1}{2^{k+1}}\right)^{(2t-2\pi ij)/\ln 2}\right]\Bigg).$$

Hence, for any arbitrary $d > 0$,

$$\begin{aligned}
P(t,z) = &-\frac{z(1-e^t)e^{2t}}{4(1-2e^{2t})\ln 2} \\
&\times \left[4e^t \sum_{k=0}^{\infty} \frac{e^{tk}(1-e^t)}{2^{k+1}+1} - 3e^t - 1\right] \\
&- \frac{1}{2}z^2 e^{2t} - \frac{z^{2+2t/\ln 2}\Gamma(-2t/\ln 2)e^{3t}(e^t-1)}{(1-2e^{2t})\ln 2} \\
&\times \Biggl\{\frac{1}{2}(e^{2t}-e^t-2) \\
&- (e^t-1)\sum_{k=0}^{\infty} e^{tk}\left[1-\left(1+\frac{1}{2^{k+1}}\right)^{2t/\ln 2}\right]\Biggr\} \\
&+ \kappa(t,\lg z) + O(z^{-d}),
\end{aligned}$$

where

$$\begin{aligned}
\kappa(t,\lg z) = &(e^t-1) \\
&\times \sum_{j\neq 0}\Biggl(\frac{2\pi i j z\Gamma(2\pi i j/\ln 2)e^{-2\pi i j\lg z}}{(1-2e^{2t})\ln^2 2} \\
&\times \Biggl(\frac{e^{2t}(2-e^t)}{1-2e^t} \\
&+ e^{3t}(e^t-1)\sum_{k=0}^{\infty} e^{tk}\left[1-\left(1+\frac{1}{2^{k+1}}\right)^{-1-2\pi i j/\ln 2}\right]\Biggr) \\
&- \frac{z^{2+2t/\ln 2}e^{-2\pi i j\lg z}\Gamma((-2t+2\pi i j)/\ln 2)}{(1-2e^{2t})\ln 2} \\
&\times \Biggl(\frac{e^{2t}-e^t-2}{2} \\
&- e^{3t}(e^t-1)\sum_{k=0}^{\infty} e^{tk}\left[1-\left(1+\frac{1}{2^{k+1}}\right)^{(2t-2\pi i j)/\ln 2}\right]\Biggr)\Biggr).
\end{aligned}$$

Recall that $P(t,z) = e^{-z}\Phi(t,z) - e^{2t}\frac{z^2}{2}$ to get

$$\mathbf{E}\left[\binom{N(z)}{2} e^{\Delta_{N(z)}t}\right]$$

$$= \frac{z(e^t-1)e^{2t}}{4(1-2e^{2t})\ln 2}\left[4e^t\sum_{k=0}^{\infty}\frac{e^{tk}(1-e^t)}{2^{k+1}+1}-3e^t-1\right]$$
$$-\frac{z^{2+2t/\ln 2}\Gamma(-2t/\ln 2)e^{3t}(e^t-1)}{(1-2e^{2t})\ln 2}$$
$$\times\left\{\frac{1}{2}(e^{2t}-e^t-2)-(e^t-1)\sum_{k=0}^{\infty}e^{tk}\left[1-\left(1+\frac{1}{2^{k+1}}\right)^{2t/\ln 2}\right]\right\}$$
$$+\kappa(t,\lg z)+O(z^{-d}).$$

We next check the conditions of the depoissonization lemma [7]. For the reader's convenience, we state this lemma.

Lemma 2 ([7]). *Let* $\{a_n\}_{n=0}^{\infty}$ *be a sequence of real numbers. Suppose that the poissonized function* $\mathcal{P}(z)=\sum_{j=0}^{\infty}\frac{a_j z^j}{j!}e^{-z}$ *exists and can be analytically continued as an entire function of complex* $z$. *Fix* $\theta\in(0,\pi/2)$ *and let* $S_\theta$ *be the cone* $\{z: |\arg z|\le\theta\}$. *Suppose that there exist positive constants* $\alpha<1$, $\beta_1$, $\beta_2$, $c$ *and* $z_0$ *such that the following conditions hold simultaneously*:

(i) *For all* $z\in S_\theta$ *with* $|z|\ge z_0$,
$$|\mathcal{P}(z)|\le\beta_1|z|^c.$$

(ii) *For all* $z\notin S_\theta$ *with* $|z|\ge z_0$,
$$|\mathcal{P}(z)e^z|\le\beta_2|z|^c e^{\alpha|z|}.$$

*Then for large* $n$,
$$a_n=\mathcal{P}(n)+O(n^{c-1/2}\ln n).$$

For the rest of the paper we restrict $t$ to the interval $[-0.1,0.1]$. For any fixed $t$ in this range, we can think of $\phi_n(t)$ as a sequence of $n$. It is clear that for this range of $t$, the poissonized function $P(t,z)$ is $O(z^{2+2t/\ln 2})$ inside the cone $S_\theta$, and is $O(|z|^{2+2t/\ln 2}e^{|z|\cos\theta})$ outside that cone.

By Lemma 2, $\mathbf{E}[\binom{n}{2}e^{\Delta_n t}]$ is the same as the poissonized version, with $n$ replacing $z$. When we divide by $\binom{n}{2}$, we get a depoissonization correction error of the order $O(n^{2t/\ln 2-1/2}\lg n)$ in the moment generating function. The Mellin inversion gives an error term of $O(n^{-d})$ for arbitrary $d>0$. Thus, combined, the two errors cannot be brought below $O(n^{2t/\ln 2-1/2}\lg n)$.

Therefore, we have

$$\begin{aligned}\mathbf{E}[e^{\Delta_n t}] &= \frac{2n^{2t/\ln 2}\Gamma(-2t/\ln 2)e^{3t}(1-e^t)}{(1-2e^{2t})\ln 2} \\ &\quad\times\Bigg(\frac{1}{2}(e^{2t}-e^t-2) \\ &\qquad - (e^t-1)\sum_{k=0}^{\infty} e^{tk}\bigg[1-\Big(1+\frac{1}{2^{k+1}}\Big)^{2t/\ln 2}\bigg]\Bigg) \\ &\quad + \frac{2}{n^2}u(t,\lg n) + O(n^{2t/\ln 2-1/2}\lg n),\end{aligned} \tag{5}$$

where

$$\begin{aligned}u(t,\lg n) = \sum_{j\neq 0} &-\frac{n^{2+2t/\ln 2}e^{-2\pi ij\lg n}(e^t-1)}{(1-2e^{2t})\ln 2}\Gamma\Big(\frac{-2t+2\pi ij}{\ln 2}\Big) \\ &\times\Bigg(\frac{e^{2t}-e^t-2}{2} \\ &\qquad - e^{3t}(e^t-1)\sum_{k=0}^{\infty} e^{tk}\bigg[1-\Big(1+\frac{1}{2^{k+1}}\Big)^{(2t-2\pi ij)/\ln 2}\bigg]\Bigg).\end{aligned}$$

[Note that the rest of the terms of $\kappa(t,\lg z)$ are subsumed by the $O$ error term.]

3.3. *The mean.* It may be possible to determine the moments of $\Delta_n$ directly from (5). However, the error requires subtle handling owing to the presence of $O$ terms, which cannot be differentiated without some regularity conditions. It is about the same effort to bypass these regularity conditions, and work indirectly from the Mellin transforms of the moments.

The derivative $\partial^k P^*(t,s)/\partial t^k$, evaluated at $t=0$, yields (a shifted version of) the Mellin transform of the expectation of the poissonized quantity $\binom{N(z)}{2}\Delta^k_{N(z)}$. We start with the mean ($k=1$). Let $e^{-z}A(z)=\mathbf{E}[\binom{N(z)}{2}\Delta_{N(z)}]$, and $B(z)=e^{-z}A(z)-z^2$. Taking the first derivative of $P^*(t,s)$ (cf. Lemma 1) and evaluating it at $t=0$, we get

$$B^*(s) = -\frac{2^{s+1}\Gamma(s+2)}{(1-2^{s+1})(1-2^{s+2})},$$

well defined for $\langle -3,-2\rangle$. The inverse Mellin transform is

$$B(z) = \frac{1}{2\pi i}\int_{-5/2-i\infty}^{-5/2+i\infty} -\frac{2^{s+1}\Gamma(s+2)}{(1-2^{s+1})(1-2^{s+2})}z^{-s}\,ds.$$

The transform $B^*(s)$ has a double pole at $s=-2$, which will provide the asymptotically dominant term (the major part of the complex integration), and simple poles at $s = -1 + \frac{2\pi ki}{\ln 2}$, for $k = 0, \pm 1, \dots$, at $s = -2 + \frac{2\pi ji}{\ln 2}$, for $j = \pm 1, \dots$, and at $s = -3, -4, \dots$. The corresponding graph with the poles of $B^*(s)$ and the area where the function is well defined is depicted in Figure 3.

We use again the method of closing the box, as was done for the moment generating function, to argue that the inverse Mellin transform of $B^*(s)$ is given by the sum of the residues of poles lying to the right of the vertical line $\Re s = -\frac{5}{2}$ and a correction term of the order $O(z^{-d})$.

We have the following residues:

$$\underset{s=-2}{Res}[B^*(s)z^{-s}] = -z^2 \lg z + z^2 \frac{3\ln 2 - 2\gamma}{2\ln 2},$$

$$\underset{s=-1}{Res}[B^*(s)z^{-s}] = -\frac{z}{\ln 2},$$

$$\sum_{k\neq 0} \underset{s=-1+2\pi ki/\ln 2}{Res}[B^*(s)z^{-s}] = -z\eta_1(\lg z)$$

and

$$\sum_{j\neq 0} \underset{s=-2+2\pi ji/\ln 2}{Res}[B^*(s)z^{-s}] = z^2\eta_2(\lg z),$$

where $\gamma = 0.577215\dots$ is Euler's constant,

$$\eta_1(u) := \frac{1}{\ln 2}\sum_{k\neq 0}\Gamma\left(1 + \frac{2\pi ik}{\ln 2}\right)e^{-2\pi iku}$$

and

$$\eta_2(u) := \frac{1}{\ln 2}\sum_{j\neq 0}\Gamma\left(\frac{2\pi ij}{\ln 2}\right)e^{-2\pi iju}. \tag{6}$$

Remark. The series in the functions $\eta_1$ and $\eta_2$ are absolutely convergent, in view of the fast asymptotic decay of the gamma function $\Gamma(x+iy)$, for fixed $x$ and increasing $y$—it is well known that $\Gamma(x+iy)$ is $O(y^{x-1/2}e^{-\pi y/2})$. For instance, for some positive constant $K$, we have

$$\begin{aligned}
|\eta_2(u)| &\le \frac{1}{\ln 2}\sum_{j\neq 0}\left|\Gamma\left(\frac{2\pi ij}{\ln 2}\right)\right| \\
&= O(1) + \frac{K}{\ln 2}\sum_{j=1}^{\infty}\left(\frac{2\pi j}{\ln 2}\right)^{-1/2} e^{-\pi^2 j/\ln 2} \\
&< O(1) + \frac{K}{\sqrt{2\pi \ln 2}}\sum_{j=1}^{\infty} e^{-\pi^2 j/\ln 2} \\
&= O(1) + \frac{K}{(1 - e^{-\pi^2/\ln 2})\sqrt{2\pi\ln 2}}.
\end{aligned}$$

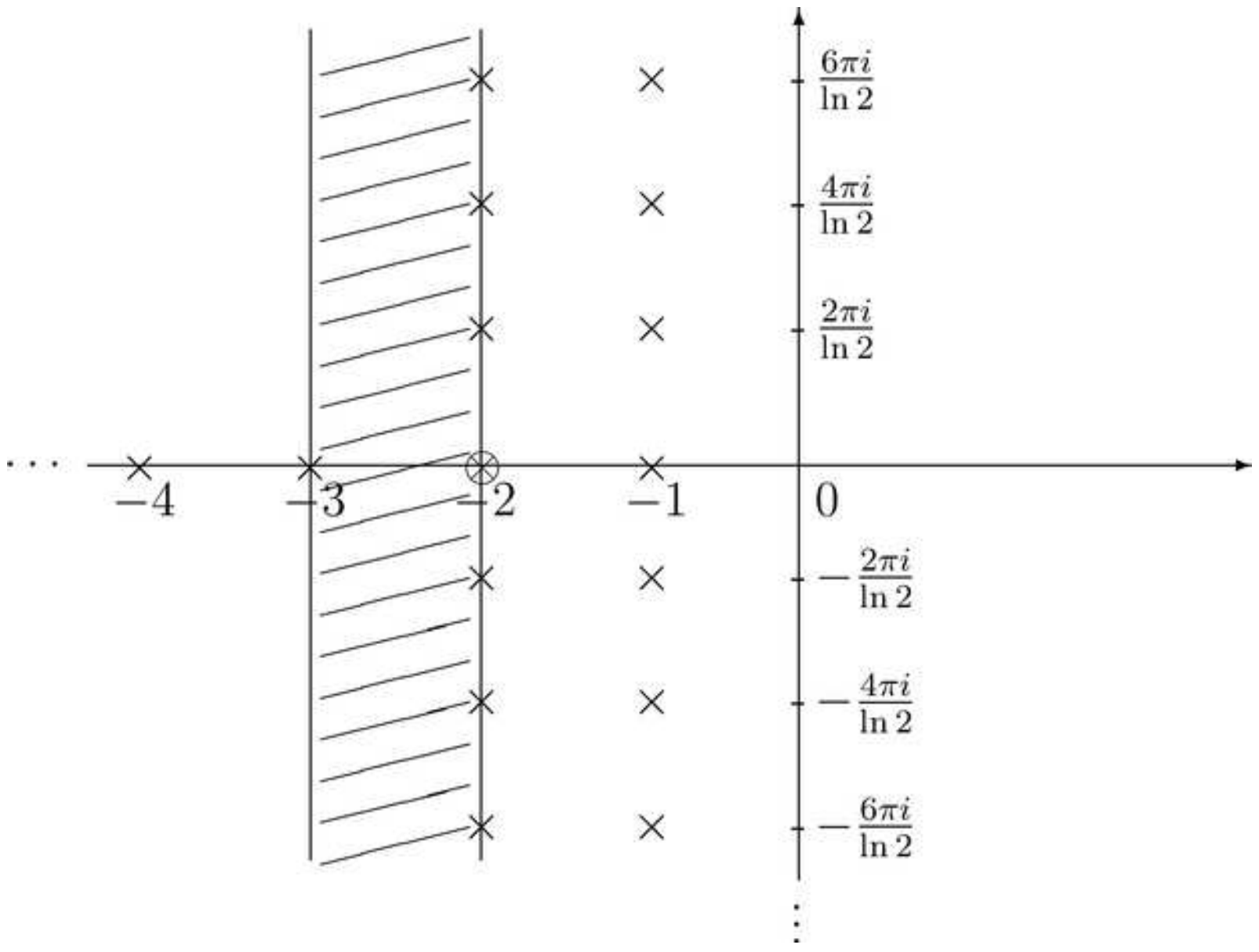

FIG. 3. *The poles of the Mellin transform $B^*(s)$; $\times$ indicates a simple pole and $\otimes$ indicates a double pole. The shaded area is the domain where $B^*(s)$ is well defined.*

Hence, for any arbitrary $d > 0$,

$$B(z) = z^2 \lg z - z^2 \eta_2(\lg z) - z^2 \frac{3\ln 2 - 2\gamma}{2\ln 2} + z\left[\frac{1}{\ln 2} + \eta_1(\lg z)\right] + O(z^{-d}).$$

The two functions $\eta_1$ and $\eta_2$ are oscillating functions of a truly ignorable magnitude. For all $z$, the magnitude of $\eta_1$ is bounded by $0.1426024772 \times 10^{-4}$ and the magnitude of $\eta_2$ by $0.1573158429 \times 10^{-5}$. The oscillations occur in the lower-order terms, anyway, and do not make $B(z)$ behave in any unstable way. This means that, as a function of $z$, the oscillations are just "waves" of a very small amplitude, fluctuating around a "steady" component.

Recall that $B(z) = e^{-z}A(z) - z^2$, and that $A(z) = \sum_{n=0}^{\infty} \binom{n}{2} \frac{\mathbf{E}[\Delta_n]}{n!} z^n$. So,

$$\mathbf{E}\left[\binom{N(z)}{2} \Delta_{N(z)}\right] = z^2 \lg z - z^2 \eta_2(\lg z) - z^2 \frac{\ln 2 - 2\gamma}{2\ln 2} + z\left[\frac{1}{\ln 2} + \eta_1(\lg z)\right] + O(z^{-d}).$$

The right-hand side in the latter relation is $O(z^{2+\varepsilon})$, for any $\varepsilon > 0$, and as in the case of the moment generating function, one can check that the conditions of the depoissonization lemma [7] hold. So, in depoissonized form,

the result is

$$\mathbf{E}\left[\binom{n}{2}\Delta_n\right] = n^2 \lg n - n^2\eta_2(\lg n) - n^2\frac{\ln 2 - 2\gamma}{2\ln 2} + n\left[\frac{1}{\ln 2} + \eta_1(\lg n)\right] + O(\max\{n^{-d}, n^{3/2+\varepsilon}\ln n\}),$$

for arbitrary $\varepsilon > 0$. As $d > 0$ is also arbitrary in the error term, we can take the error to be $O(n^{3/2+\varepsilon})$. When we divide by $\binom{n}{2}$, we get $O(n^{-1/2+\varepsilon})$ error term. For ease of exposition, we take $\varepsilon = 0.0001$ and present the error in the order $O(n^{-0.4999})$.

PROPOSITION 1. *In a trie of $n$ random keys following the unbiased Bernoulli model, the distance $\Delta_n$ between two randomly selected keys has the mean value*

$$\mathbf{E}[\Delta_n] = 2\lg n + \eta(\lg n) - \frac{\ln 2 - 2\gamma}{\ln 2} + O\left(\frac{1}{n^{0.4999}}\right),$$

*where $\eta(\cdot) = -2\eta_2(\cdot)$, with $\eta_2(\cdot)$ being the small oscillating function given in (6).*

3.4. *The variance.* One can continue pumping moments by taking derivatives of $P^*(t,s)$. The derivation of higher moments follows the same general principles that were used for the mean: fixed-population recurrence, poissonization, asymptotic solution by the Mellin transform and its inverse, then depoissonization. However, the task becomes too daunting, even for the second moment. This combinatorial explosion phenomenon is common folklore in random graphs. Let

$$L(z) := \mathbf{E}\left[\binom{N(z)}{2}\Delta^2_{N(z)}\right] - 2z^2.$$

Hence,

$$\begin{aligned} L^*(s) &= \left.\frac{\partial^2 P^*(t,s)}{\partial t^2}\right|_{t=0} \\ &= \frac{2^{s+1}\Gamma(s+2)}{(1-2^{s+1})(1-2^{s+2})} \\ &\quad \times \left(-\frac{2^{s+4}}{1-2^{s+2}} + \frac{1-7\cdot 2^{s+2}}{2^{s+2}} - 2\sum_{k=0}^{\infty}\left[1 - \frac{1}{(1+1/2^{k+1})^{s+2}}\right]\right), \end{aligned}$$

well defined in $\langle -3, -2\rangle$.

The inverse Mellin transform is

$$L(z) = \frac{1}{2\pi i}\int_{-5/2-i\infty}^{-5/2+i\infty} L^*(s) z^{-s}\, ds,$$

which again we handle by shifting the line of integration to the right and account for the shift by the residue of the poles introduced. The transform $L^*(s)$ has a triple pole at $s=-2$, which provides the asymptotically dominant term, double poles at $s=-2+\frac{2\pi ji}{\ln 2}$, for $j=\pm 1,\dots$, and simple poles at $s=-1+\frac{2\pi ki}{\ln 2}$, for $k=0,\pm 1,\dots$ and at $s=-3,-4,\dots$.

We have the following residue contributions:

$$\begin{aligned}
\underset{s=-2}{Res}[L^*(s)z^{-s}] &= -2z^2\lg^2 z - 2z^2\lg z\left(\frac{2\gamma}{\ln 2}-1\right) \\
&\quad - \frac{z^2}{\ln^2 2}\left[\frac{1}{3}\pi^2 + 2\gamma^2 - 2\gamma\ln 2 + \frac{5}{3}\ln^2 2\right. \\
&\quad \left. - 2\ln 2\sum_{k=0}^{\infty}\ln\left(1+\frac{1}{2^{k+1}}\right)\right], \\
\underset{s=-1}{Res}[L^*(s)z^{-s}] &= \frac{z}{\ln 2}\left(\frac{3}{2} - 2\sum_{k=0}^{\infty}\frac{1}{1+2^{k+1}}\right),
\end{aligned}$$

$$\sum_{j\neq 0}\underset{s=-2+2\pi ji/\ln 2}{Res}[L^*(s)z^{-s}] = -z^2\xi_1(\lg z),$$

$$\sum_{k\neq 0}\underset{s=-1+2\pi ki/\ln 2}{Res}[L^*(s)z^{-s}] = z\xi_2(\lg z),$$

where $\xi_1(\cdot)$ an $\xi_2(\cdot)$ are oscillating function given by

$$\begin{aligned}
\xi_1(u) &:= \frac{1}{\ln 2}\sum_{j\neq 0}\Gamma\left(\frac{2\pi ij}{\ln 2}\right)e^{-2\pi iju} \\
&\quad \times\left(10 - 2\sum_{k=0}^{\infty}\left[1 - \frac{1}{(1+1/2^{k+1})^{2\pi ij/\ln 2}}\right]\right. \\
&\quad \left. + \frac{4\psi(2\pi ij/\ln 2)}{\ln 2} - 4u\right), \\
\xi_2(u) &:= \frac{1}{\ln 2}\sum_{k\neq 0}\Gamma\left(1+\frac{2\pi ik}{\ln 2}\right)e^{-2\pi iku} \\
&\quad \times\left(\frac{3}{2} - 2\sum_{j=0}^{\infty}\left[1 - \frac{1}{(1+1/2^{j+1})^{1+2\pi ik/\ln 2}}\right]\right),
\end{aligned}$$

and $\psi(\cdot)$ is the digamma function.

Hence, for any arbitrary $\theta > 0$,

$$L(z) = 2z^2\lg^2 z + 2z^2\lg z\left(\frac{2\gamma}{\ln 2}-1\right)$$

$$+ \frac{z^2}{\ln^2 2}\left[\frac{1}{3}\pi^2 + 2\gamma^2 - 2\gamma \ln 2 + \frac{5}{3}\ln^2 2 - 2\ln 2 \sum_{k=0}^{\infty} \ln\left(1 + \frac{1}{2^{k+1}}\right)\right]$$

$$- \frac{z}{\ln 2}\left(\frac{3}{2} - 2\sum_{k=0}^{\infty} \frac{1}{1+2^{k+1}}\right) + z^2\xi_1(\lg z) - z\xi_2(\lg z) + O(z^{-\theta}).$$

Recall that $L(z)$ is a shifted poissonized version of the problem, with a number of keys $N(z)$ that follows a Poisson distribution with parameter $z$. Therefore,

$$\begin{aligned} \mathbf{E}\left[\binom{N(z)}{2}\Delta^2_{N(z)}\right] &= 2z^2 \lg^2 z + 2z^2 \lg z\left(\frac{2\gamma}{\ln 2} - 1\right) \\ &\quad + \frac{z^2}{\ln^2 2}\left[\frac{1}{3}\pi^2 + 2\gamma^2 - 2\gamma\ln 2 + \frac{11}{3}\ln^2 2 \right. \\ &\qquad\qquad \left. - 2\ln 2\sum_{k=0}^{\infty}\ln\left(1+\frac{1}{2^{k+1}}\right)\right] \\ &\quad - \frac{z}{\ln 2}\left(\frac{3}{2} - 2\sum_{k=0}^{\infty}\frac{1}{1+2^{k+1}}\right) \\ &\quad + z^2\xi_1(\lg z) - z\xi_2(\lg z) + O(z^{-\theta}). \end{aligned} \tag{7}$$

The depoissonization of (7) by the depoissonization lemma [7] finally gives us the second moment for the fixed population in the same form with $n$ replacing $z$; the depoissonization introduces a small error of the order $O(n^{-0.4999})$; the genesis of this particular representation of error is discussed in the paragraphs preceding Proposition 1. The entire validity of depoissonization bears on arguments similar to those discussed in some more detail prior to Proposition 1. The variance follows by subtracting the square of the mean (cf. Proposition 1). Interestingly, the $\lg^2 n$ and $\lg n$ terms both disappear, leaving behind only a $O(1)$ function, consisting of a steady component and negligible oscillations.

THEOREM 1. *In a trie of $n$ random keys following the unbiased Bernoulli model, the distance $\Delta_n$ between two randomly selected keys has the variance*

$$\mathbf{Var}[\Delta_n] = \frac{2\pi^2 + 19\ln^2 2 - 12\ln 2\sum_{k=0}^{\infty}\ln(1+1/2^{k+1})}{3\ln^2 2} + \xi(\lg n) + O\left(\frac{1}{n^{0.4999}}\right),$$

*where $\xi(\cdot)$ is an oscillating function given by*

$$\xi(u) = \frac{16\ln 2 + 8\gamma}{\ln 2}\eta_2(u) - 4\eta_2^2(u) - 4\xi_3(u),$$

*with* $\eta_2(\cdot)$ *given in* (6) *and* $\xi_3(\cdot)$ *given by*

$$\xi_3(u) := \frac{1}{\ln 2}\sum_{j\neq 0}\Gamma\left(\frac{2\pi ij}{\ln 2}\right)e^{-2\pi iju} \times\left(\sum_{k=0}^{\infty}\left[1-\frac{1}{(1+1/2^{k+1})^{2\pi ij/\ln 2}}\right]-\frac{2\psi(2\pi ij/\ln 2)}{\ln 2}\right).$$

The function $\xi_3$ in Theorem 1 is absolutely convergent, by arguments that are essentially the same as those given for the function $\eta_1$ and $\eta_2$ that appear in the mean; see the remark following (6). In fact, $\xi$ offers only negligible oscillations; its magnitude is bounded by $0.5654042648\times 10^{-4}$. The arguments in [16] may suggest a way to simplify $\xi$. In view of the sharp $O(1)$ variance concentration of $\Delta_n$, a concentration law is an immediate corollary (by Chebyshev's inequality).

Corollary 1. *As* $n\to\infty$,

$$\frac{\Delta_n}{\lg n}\xrightarrow{\mathcal{P}} 2.$$

3.5. *The oscillatory distribution.* The presence of the oscillating function $\kappa(t,\lg n)$ in the moment generating function of $\Delta_n$ indicates that a limit distribution does not exist. In this section we make this argument more explicit. We have determined in Section 3.3 that the mean has $2\lg n$ leading asymptotic term. We have also seen in Section 3.4 that the variance is $O(1)$. One then looks into the behavior of the integer random variable of $\Delta_n^* = \Delta_n - \lfloor 2\lg n\rfloor$, if there is a chance for convergence. The main result of this paper answers this question in the negative.

Theorem 2. *Let* $\Delta_n$ *be the distance between two randomly selected keys in a trie of* $n$ *random keys following the unbiased Bernoulli model. Then* $\Delta_n - \lfloor 2\lg n\rfloor$ *has an oscillating moment generating function, and, consequently, does not converge in distribution to any random variable.*

Proof. One finds from (5) that $\Delta_n^*$ has the moment generating function

$$\mathbf{E}[e^{\Delta_n^* t}] = (G(t)+H_n(t))e^{t(2\lg n-\lfloor 2\lg n\rfloor)} + O\left(\frac{\lg n}{\sqrt{n}}e^{t(2\lg n-\lfloor 2\lg n\rfloor)}\right), \tag{8}$$

where

$$G(t) = \frac{2\Gamma(-2t/\ln 2)e^{3t}(1-e^t)}{(1-2e^{2t})\ln 2} \times\left(\frac{1}{2}(e^{2t}-e^t-2)-(e^t-1)\sum_{k=0}^{\infty}e^{tk}\left[1-\left(1+\frac{1}{2^{k+1}}\right)^{2t/\ln 2}\right]\right)$$

and

$$H_n(t) = \sum_{j\neq 0} -\frac{2e^{-2\pi i j \lg n}(e^t-1)}{(1-2e^{2t})\ln 2}\Gamma\left(\frac{-2t+2\pi i j}{\ln 2}\right) \times \left(\frac{e^{2t}-e^t-2}{2} - e^{3t}(e^t-1)\sum_{k=0}^{\infty} e^{tk}\left[1-\left(1+\frac{1}{2^{k+1}}\right)^{(2t-2\pi i j)/\ln 2}\right]\right).$$

Let $\{x\}$ denote the fractional part of $x$, that is,

$$\{x\} = x - \lfloor x \rfloor.$$

Using this notation in (8), we find

$$\mathbf{E}[e^{\Delta_n^* t}] = (G(t) + H_n(t))e^{\{2\lg n\}t} + O\left(\frac{\lg n}{\sqrt{n}} e^{\{2\lg n\}t}\right),$$

for the range $|t| < \frac{1}{2}\ln 2$. For the rest of the proof, we restrict $t$ to the interval $[-0.1, 0.1]$. For any fixed $t$ from this range, the function $H_n(t)$ provides small ignorable oscillations around $G(t)$, and, of course, the $O(n^{-1/2}\lg n e^{\{2\lg n\}t})$ term can be made arbitrarily small. It is well known that the function $\{2\lg n\}$ is dense in the interval $[0, 1)$; see, for example, [10]. The term $\{2\lg n\}$ is genuinely oscillating. Also, for $t > 0$,

$$\inf_{n\geq 0} G(t)e^{\{2\lg n\}t} = G(t) \quad \text{and} \quad \sup_{n\geq 0} G(t)e^{\{2\lg n\}t} = e^t G(t)$$

and, for $t < 0$,

$$\inf_{n\geq 0} G(t)e^{\{2\lg n\}t} = e^t G(t) \quad \text{and} \quad \sup_{n\geq 0} G(t)e^{\{2\lg n\}t} = G(t).$$

Then, for large $n$, the moment generating function $\mathbf{E}[e^{\Delta_n^* t}]$ itself comes infinitely often very close to the two functions $G(t)$ and $G(t)e^t$. These two "envelopes" are distinct for all values of $t$ with the exception of that of $t = 0$, where both values are equal to $G(0) = 1$. Namely,

$$|(G(t) + H_n(t))e^{\{2\lg n\}t}| \geq (G(t) - |H_n(t)|)e^{\{2\lg n\}t},$$

and for nonnegative $t$ in the range of interest, we have

$$|H_n(t)| \leq \frac{2|1-e^t|}{|1-2e^{2t}|\ln 2} \times \left(\left|\frac{e^{2t}-e^t-2}{2}\right| + e^{3t}|e^t-1|\sum_{k=0}^{\infty} e^{tk}\left[\left(1+\frac{1}{2^{k+1}}\right)^{2t/\ln 2} - 1\right]\right)$$

$$
\begin{aligned}
&\times \sum_{j\neq 0}\left|\Gamma\left(\frac{-2t+2\pi i j}{\ln 2}\right)\right| \\
&\le \frac{2|1-e^t|}{|1-2e^{2t}|\ln 2} \\
&\quad \times \left(\left|\frac{e^{2t}-e^t-2}{2}\right| + e^{3t}|e^t-1|\sum_{k=0}^{\infty} e^{tk}\left[\left(1+\frac{1}{2^{k+1}}\right)-1\right]\right) \\
&\quad \times \sum_{j\neq 0}\left|\Gamma\left(\frac{-2t+2\pi i j}{\ln 2}\right)\right| \\
&\le \frac{2|1-e^t|}{|1-2e^{2t}|\ln 2}\left(\left|\frac{e^{2t}-e^t-2}{2}\right| + \frac{e^{3t}|1-e^t|}{|2-e^t|}\right)\sum_{j\neq 0}\left|\Gamma\left(\frac{-2t+2\pi i j}{\ln 2}\right)\right| \\
&\le 3\times 10^{-4}.
\end{aligned}
$$

The latter numerical bound follows from maximizing the part outside the sum by standard calculus arguments, and working out a numerical value for the sum as was done for the mean [see the remarks following (6)]. Similarly, for negative $t$ in the range of interest,

$$
\begin{aligned}
|H_n(t)| &\le \frac{2|1-e^t|}{|1-2e^{2t}|\ln 2} \\
&\quad \times \left(\left|\frac{e^{2t}-e^t-2}{2}\right| + e^{3t}|e^t-1|\sum_{k=0}^{\infty} e^{tk}\left[1-\left(1+\frac{1}{2^{k+1}}\right)^{2t/\ln 2}\right]\right) \\
&\quad \times \sum_{j\neq 0}\left|\Gamma\left(\frac{-2t+2\pi i j}{\ln 2}\right)\right| \\
&\le \frac{2|1-e^t|}{|1-2e^{2t}|\ln 2}\left(\left|\frac{e^{2t}-e^t-2}{2}\right| + e^{3t}|e^t-1|\sum_{k=0}^{\infty} e^{tk}\right) \\
&\quad \times \sum_{j\neq 0}\left|\Gamma\left(\frac{-2t+2\pi i j}{\ln 2}\right)\right| \\
&\le \frac{2|1-e^t|}{|1-2e^{2t}|\ln 2}\left(\left|\frac{e^{2t}-e^t-2}{2}\right| + e^{3t}\right)\sum_{j\neq 0}\left|\Gamma\left(\frac{-2t+2\pi i j}{\ln 2}\right)\right| \\
&\le 3\times 10^{-4}.
\end{aligned}
$$

The numerical bound is found in a manner similar to that in the case of $t\ge 0$.

In either case, $3 \times 10^{-4}$ is an upper bound. Hence, for $t > t_0$ [where $0 < t_0 < 0.1$ is the point where $(G(t) - 0.0003)e^t$ and $G(t) + 0.0003$ intersect],

$$\limsup_{n\to\infty} \mathbf{E}[e^{\Delta_n^* t}] \ge \limsup_{n\to\infty} (G(t) - 0.0003)e^{\{2\lg n\}t} = (G(t) - 0.0003)e^t.$$

Similarly,

$$\liminf_{n\to\infty} \mathbf{E}[e^{\Delta_n^* t}] \le G(t) + 0.0003.$$

The two bounds are different, as $G(t) \ge 0.98$, uniformly in $t$. For instance, letting $t = \frac{1}{10}$ gives

$$\liminf_{n\to\infty} \mathbf{E}[e^{(1/10)\Delta_n^*}] \le 1.148 < 1.267 \le \limsup_{n\to\infty} \mathbf{E}[e^{(1/10)\Delta_n^*}]. \qquad \square$$

In principle, one can invert a moment generating function. For integer $r$, we have

$$\mathbf{Pr}(\Delta_n^* = r) = \frac{1}{2\pi} \int_0^{2\pi} e^{-iur} \mathbf{E}[e^{iu\Delta_n^*}]\, du.$$

The integrals may not be easy to work through. For example, for $r = 0$, we have

$$\mathbf{Pr}(\Delta_n^* = 0) = \frac{1}{2\pi} \int_0^{2\pi} G(iu)e^{\{2\lg n\}iu}\, du + \frac{1}{2\pi} \int_0^{2\pi} H_n(iu)e^{\{2\lg n\}iu}\, du + O\left(\frac{\lg n}{\sqrt{n}}\right).$$

The functions $G(\cdot)$ and $H_n(\cdot)$ are complicated, but it is clear that $\mathbf{Pr}(\Delta_n^* = 0)$ is oscillating.

**4. Discussion.** Tries offer many applications in different areas, such as computer science, telecommunications and computational biology. Applications include a variety of pattern matching algorithms on words, as well as a model for analyzing the periodicities and autocorrelation between the substrings of a string. This is an area of interest to the study of DNA sequences (digital data on an alphabet of four letters) in computational biology, when the sequences are almost random.

By purely analytic methods, we investigated distances between pairs of keys in a random trie on $n$ keys. The average value is $2\lg n$ (modulated by some small oscillations), the variance is $O(1)$, and the centered distance does not have a limit distribution, but rather oscillates between two extremal values. The nonexistence of limits in the context of digital trees has been noted before in [3], who studied the height of random incomplete tries by a mix of probabilistic and analytic methods. The probabilistic behavior of other related parameters can be obtained from our results. For example, it is a corollary to our results that the Wiener index

$$W_n = \sum D_{ij},$$

the sum of all pairs of distances ($D_{ij}$ is the distance between the $i$th and $j$th keys), has an average value equal to

$$\mathbf{E}[W_n] = n^2 \lg n - \left(\frac{\ln 2 - 2\gamma}{2\ln 2} + \eta_2(\lg n)\right) n^2 + O(n^{1.5001}).$$

To put this result in appropriate perspective, compare with [14].

A trie is said to be *complete* when it is fully grown at each level (except possibly the last one). All the levels in a complete trie are filled with internal nodes, except the last two (or one) where the keys reside. Thus, in such a complete trie, the keys are at distance $\lfloor \lg n \rfloor$ or $1 + \lfloor \lg n \rfloor$. This balance is usually a desired feature for the fast retrieval of data on average. In the complete trie, if pairs of leaves are chosen at random, both will always be at depth $\lg n + O(1)$. It can be shown by a recursive argument that the distance between a pair of keys has an asymptotic average $2 \lg n$ and variance $O(1)$. Thus, our results confirm that the trie is an excellent choice for random digital data; it tends to naturally balance itself, behaving rather closely as an ideally balanced complete trie.

It appears to us that the results can be extended without introducing essential difficulty to the case of data on alphabets larger than binary, which can be useful for DNA studies. The problem that can possibly be of different structure, giving an interesting offshoot, is that where the Bernoulli model is biased. In the binary case, for example, one would consider an ergodic source emitting bits of 1's with probability $p \neq \frac{1}{2}$, and 0's with probability $q = 1 - p \neq \frac{1}{2}$. Several symmetries employed in the balanced Bernoulli case cease to exist. Under a biased Bernoulli model, the functional equations obtained, both under a fixed-population model and a poissonized data model, will be markedly different. We leave this as a topic for future research.

## APPENDIX A

**A Mellin viewpoint of the random depth.** There are several known results for the random depth derived by probabilistic methods [see [1] for exact analysis]. See also the related studies of asymptotics in [6, 8, 17]; these results are surveyed in [11]. However, we need the Mellin transform viewpoint of these results to plug in the functional equations arising for $\Delta_n$.

We start from

$$\delta_n | L_n = \begin{cases} \delta_{L_n} + 1, & \text{with probability } \dfrac{L_n}{n}, \\ \tilde{\delta}_{R_n} + 1, & \text{with probability } \dfrac{R_n}{n}. \end{cases}$$

As in the formulation of the functional equation (2) for $\Delta_n$, we use the symmetry of the left and the right subtrees and the fact that $L_n$ is distributed like $\mathrm{Bin}(n, \frac{1}{2})$. After multiplying by $z^n$ and summing over all possible values of $n$, we get

$$\Psi(t,z) = 2e^t e^{z/2}\Psi\left(t, \frac{z}{2}\right) + z(1-e^t),$$

where $\Psi(t,z) := \sum_{j=0}^{\infty} z^j \frac{j\phi_{\delta_j}(t)}{j!}$. Multiplying further by $e^{-z}$ to poissonize, we get

$$e^{-z}\Psi(t,z) = 2e^t e^{-z/2}\Psi\left(t, \frac{z}{2}\right) + ze^{-z}(1-e^t).$$

We then introduce a shift of $z$ (to ensure the existence of the Mellin transform) to get

$$e^{-z}\Psi(t,z) - z = 2e^t\left[e^{-z/2}\Psi\left(t, \frac{z}{2}\right) - \frac{z}{2}\right] + ze^t - z + ze^{-z}(1-e^t),$$

that is,

$$Q(t,z) = 2e^t Q\left(t, \frac{z}{2}\right) - z(1-e^{-z})(1-e^t), \tag{9}$$

where $Q(t,z) := e^{-z}\Psi(t,z) - z$.

One can recover the asymptotic result of Jacquet and Régnier [6] or Pittel [17] by completing a Mellin inversion and a depoissonization program. However, our purpose is specific to the task of developing a functional equation for

$$\Theta(t,z) := Q^2(t,z),$$

as we need its transform in $P^*(t,s)$. Hence, square (9) above to get

$$\begin{aligned} Q^2(t,z) &= 4e^{2t}Q^2\left(t, \frac{z}{2}\right) \\ &\qquad - 4e^t z(1-e^{-z})(1-e^t)Q\left(t, \frac{z}{2}\right) + z^2(1-e^{-z})^2(1-e^t)^2, \end{aligned}$$

that is,

$$\begin{aligned} \Theta(t,z) &= 4e^{2t}\Theta\left(t, \frac{z}{2}\right) \\ &\qquad - 4e^t z(1-e^{-z})(1-e^t)Q\left(t, \frac{z}{2}\right) + z^2(1-e^{-z})^2(1-e^t)^2. \end{aligned}$$

Therefore,

$$\begin{aligned} \Theta^*(t,s) &= 4e^{2t}2^s\Theta^*(t,s) - \mathbf{Mellin}\left\{4e^t z(1-e^{-z})(1-e^t)Q\left(t, \frac{z}{2}\right); s\right\} \\ &\qquad + \mathbf{Mellin}\{z^2(1-e^{-z})^2(1-e^t)^2; s\}, \end{aligned}$$

which means that

$$\Theta^*(t,s) = \Big(-\mathbf{Mellin}\Big\{4e^t z(1-e^{-z})(1-e^t)Q\Big(t,\frac{z}{2}\Big);s\Big\} + \mathbf{Mellin}\{z^2(1-e^{-z})^2(1-e^t)^2;s\}\Big)\Big/(1-e^{2t}2^{s+2}).$$

## APPENDIX B

**Mellin transform of $P(t,z)$.** Taking the Mellin transform of (3), we get

$$P^*(t,s) = 2^{s+1}P^*(t,s) + \mathbf{Mellin}\Big\{F\Big(t,\frac{z}{2}\Big) - e^{2t}\frac{z^2}{4};s\Big\}.$$

That is,

$$P^*(t,s) = \frac{\mathbf{Mellin}\{F(t,z/2) - e^{2t}z^2/4;s\}}{1-2^{s+1}}. \tag{10}$$

However,

$$\begin{aligned} F\Big(t,\frac{z}{2}\Big) - e^{2t}\frac{z^2}{4} &= e^{2t}e^{-z}\Big[\frac{z}{2} + \sum_{j=2}^{\infty}\frac{j\phi_{\delta_j}(t)}{j!}\Big(\frac{z}{2}\Big)^j\Big]^2 - e^{2t}\frac{z^2}{4} \\ &= e^{2t}e^{-z}\Big[\sum_{j=2}^{\infty}\frac{j\phi_{\delta_j}(t)}{j!}\Big(\frac{z}{2}\Big)^j\Big]^2 \\ &\quad + ze^{-z}e^{2t}\sum_{j=2}^{\infty}\frac{j\phi_{\delta_j}(t)}{j!}\Big(\frac{z}{2}\Big)^j - e^{2t}\frac{z^2}{4}(1-e^{-z}). \end{aligned} \tag{11}$$

In order to get the Mellin transform of the first part of the right-hand side of (11), we first let

$$f(t,z) := e^{-z}\sum_{j=0}^{\infty}\frac{j\phi_{\delta_j}(t)}{j!}z^j.$$

Hence,

$$\begin{aligned} \Big[e^{-z/2}\sum_{j=2}^{\infty}\frac{j\phi_{\delta_j}(t)}{j!}\Big(\frac{z}{2}\Big)^j\Big]^2 &= f^2\Big(t,\frac{z}{2}\Big) + \frac{z^2}{4}e^{-z} - ze^{-z/2}f\Big(t,\frac{z}{2}\Big) \\ &= Q^2\Big(t,\frac{z}{2}\Big) + z(1-e^{-z/2})Q\Big(t,\frac{z}{2}\Big) \\ &\quad + \frac{z^2}{4}(1+e^{-z}-2e^{-z/2}), \end{aligned}$$

with $Q(t,z)=f(t,z)-z$. Therefore, for the first part on the right-hand side of (11), we have the Mellin transform

$$
\begin{aligned}
&\mathbf{Mellin}\left\{e^{2t}e^{-z}\left[\sum_{j=2}^{\infty}\frac{j\phi_{\delta_j}(t)}{j!}\left(\frac{z}{2}\right)^j\right]^2;s\right\}\\
&\quad=\mathbf{Mellin}\left\{e^{2t}Q^2\left(t,\frac{z}{2}\right);s\right\}\\
&\qquad+\mathbf{Mellin}\left\{ze^{2t}(1-e^{-z/2})Q\left(t,\frac{z}{2}\right);s\right\}\\
&\qquad+\mathbf{Mellin}\left\{\frac{z^2}{4}e^{2t}(1+e^{-z}-2e^{-z/2});s\right\}.
\end{aligned}\tag{12}
$$

Moreover, the relation

$$Q(t,z)=2e^tQ\left(t,\frac{z}{2}\right)-z(1-e^{-z})(1-e^t)$$

gives, by recursion,

$$Q\left(t,\frac{z}{2}\right)=-\sum_{k=0}^{\infty}(2e^t)^k\frac{z}{2^{k+1}}(1-e^{-z/2^{k+1}})(1-e^t).\tag{13}$$

Hence,

$$
\begin{aligned}
&\mathbf{Mellin}\left\{4e^tz(1-e^{-z})(1-e^t)Q\left(t,\frac{z}{2}\right);s\right\}\\
&\quad=-4e^t(1-e^t)^2\int_0^\infty(1-e^{-z})\sum_{k=0}^{\infty}(2e^t)^k\frac{z}{2^{k+1}}(1-e^{-z/2^{k+1}})z^s\,dz\\
&\quad=-2e^t(1-e^t)^2\sum_{k=0}^{\infty}e^{tk}\left[\int_0^\infty(1-e^{-z})z^{s+1}\,dz\right.\\
&\qquad\qquad\left.-\int_0^\infty e^{-z/2^{k+1}}z^{s+1}\,dz+\int_0^\infty e^{-z}e^{-z/2^{k+1}}z^{s+1}\,dz\right]\\
&\quad=-2e^t(1-e^t)^2\Gamma(s+2)\sum_{k=0}^{\infty}e^{tk}\left[-1-(2^{k+1})^{s+2}+\frac{1}{(1+1/2^{k+1})^{s+2}}\right].
\end{aligned}
$$

In addition,

$$
\begin{aligned}
&\mathbf{Mellin}\{z^2(1-e^{-z})^2(1-e^t)^2;s\}\\
&\quad=(1-e^t)^2\int_0^\infty z^2(1-e^{-z})^2z^{s-1}\,dz\\
&\quad=(1-e^t)^2\left[\int_0^\infty 2(1-e^{-z})z^{s+1}\,dz-\int_0^\infty(1-e^{-2z})z^{s+1}\,dz\right]\\
&\quad=(1-e^t)^2\frac{(1-2^{s+3})\Gamma(s+2)}{2^{s+2}}.
\end{aligned}
$$

Hence, we have

$$\Theta^*(t,s) = \frac{1}{1-e^{2t}2^{s+2}}\Bigg(2e^t(1-e^t)^2\Gamma(s+2) \times \sum_{k=0}^{\infty} e^{tk}\Big[-1-(2^{k+1})^{s+2} + \frac{1}{(1+1/2^{k+1})^{s+2}}\Big] + (1-e^t)^2\frac{(1-2^{s+3})\Gamma(s+2)}{2^{s+2}}\Bigg),$$

which means that the Mellin transform of the first part on the right-hand side of (12) is

$$\begin{aligned}&\mathbf{Mellin}\Big\{e^{2t}Q^2\Big(t,\frac{z}{2}\Big);s\Big\}\\ &\quad= e^{2t}2^s\Theta^*(t,s)\\ &\quad= \frac{e^{2t}2^s}{1-e^{2t}2^{s+2}}\Bigg(2e^t(1-e^t)^2\Gamma(s+2) \times \sum_{k=0}^{\infty} e^{tk}\Big[-1-(2^{k+1})^{s+2} + \frac{1}{(1+1/2^{k+1})^{s+2}}\Big] + \frac{(1-e^t)^2(1-2^{s+3})\Gamma(s+2)}{2^{s+2}}\Bigg).\end{aligned}$$

Now, using (13), we note that by similar steps we have

$$\begin{aligned}&\mathbf{Mellin}\Big\{ze^{2t}e^{-z/2}Q\Big(t,\frac{z}{2}\Big);s\Big\}\\ &\quad= \sum_{k=0}^{\infty}\frac{1}{2}e^{(k+2)t}(1-e^t) \times \int_0^{\infty} -z^2e^{-z/2}(1-e^{-z/2^{k+1}})z^{s-1}\,dz\\ &\quad= -2^{s+1}e^{2t}(1-e^t)\Gamma(s+2)\sum_{k=0}^{\infty}e^{tk}\Big[1-\frac{1}{(1+1/2^k)^{s+2}}\Big].\end{aligned}$$

Hence, the second Mellin transform on the right-hand side of (12) can be calculated as

$$\begin{aligned}&\mathbf{Mellin}\Big\{ze^{2t}(1-e^{-z/2})Q\Big(t,\frac{z}{2}\Big);s\Big\}\\ &\quad= \mathbf{Mellin}\Big\{ze^{2t}Q\Big(t,\frac{z}{2}\Big);s\Big\}\end{aligned}$$

$$
\begin{aligned}
&- \mathbf{Mellin}\Big\{ze^{2t}e^{-z/2}Q\Big(t,\frac{z}{2}\Big);s\Big\}\\
&= 2^{s+1}e^{2t}(1-e^t)\Gamma(s+2)\\
&\quad\times\sum_{k=0}^{\infty}e^{tk}\Big[2^{k(s+2)}+1-\frac{1}{(1+1/2^{k+1})^{s+2}}\Big].
\end{aligned}
$$

Finally, the third Mellin transform on the right-hand side of (12) is

$$
\begin{aligned}
&\mathbf{Mellin}\Big\{\frac{z^2}{4}e^{2t}(1+e^{-z}-2e^{-z/2});s\Big\}\\
&\quad= \mathbf{Mellin}\Big\{\frac{z^2}{4}e^{2t}[(1-e^{-z/2})-e^{-z/2}(1-e^{-z/2})];s\Big\}\\
&\quad= \Big(\frac{1}{4}-2^{s+1}\Big)e^{2t}\Gamma(s+2).
\end{aligned}
$$

Therefore, putting the three pieces together, we obtain the first part of the right-hand side of (11) as

$$
\begin{aligned}
&\mathbf{Mellin}\Big\{e^{2t}e^{-z}\Big[\sum_{j=2}^{\infty}\frac{j\phi_{\delta_j}(t)}{j!}\Big(\frac{z}{2}\Big)^j\Big]^2;s\Big\}\\
&= \frac{e^{2t}2^s}{1-e^{2t}2^{s+2}}\Big\{2e^t(1-e^t)^2\Gamma(s+2)\\
&\qquad\times\sum_{k=0}^{\infty}e^{tk}\Big[-1-(2^{k+1})^{s+2}+\frac{1}{(1+1/2^{k+1})^{s+2}}\Big]\\
&\qquad+(1-e^t)^2\frac{(1-2^{s+3})\Gamma(s+2)}{2^{s+2}}\Big\}\\
&\quad+2^{s+1}e^{2t}(1-e^t)\Gamma(s+2)\\
&\qquad\times\sum_{k=0}^{\infty}e^{tk}\Big[2^{k(s+2)}+1-\frac{1}{(1+1/2^{k+1})^{s+2}}\Big]\\
&\quad+\Big(\frac{1}{4}-2^{s+1}\Big)e^{2t}\Gamma(s+2).
\end{aligned}
\tag{14}
$$

For the second part of the right-hand side of (11), we start from

$$
\begin{aligned}
&ze^{2t}e^{-z}\sum_{j=2}^{\infty}\frac{j\phi_{\delta_j}(t)}{j!}\Big(\frac{z}{2}\Big)^j\\
&\quad= ze^{2t}e^{-z/2}f\Big(t,\frac{z}{2}\Big)-\frac{z^2}{2}e^{2t}e^{-z}\\
&\quad= ze^{2t}e^{-z/2}Q\Big(t,\frac{z}{2}\Big)+\frac{z^2}{2}e^{2t}e^{-z/2}-\frac{z^2}{2}e^{2t}e^{-z}.
\end{aligned}
$$

Hence,

$$\mathbf{Mellin}\left\{ze^{-z}e^{2t}\sum_{j=2}^{\infty}\frac{j\phi_{\delta_j}(t)}{j!}\left(\frac{z}{2}\right)^j; s\right\} = \mathbf{Mellin}\left\{ze^{2t}e^{-z/2}Q\left(t,\frac{z}{2}\right); s\right\}$$
$$+\,\mathbf{Mellin}\left\{\frac{z^2}{2}e^{2t}e^{-z/2}; s\right\}$$
$$-\,\mathbf{Mellin}\left\{\frac{z^2}{2}e^{2t}e^{-z}; s\right\}.$$

Dealing with each piece separately and utilizing results from before, we obtain the following results for the first two pieces:

$$\mathbf{Mellin}\left\{ze^{2t}e^{-z/2}Q\left(t,\frac{z}{2}\right); s\right\}$$
$$= -2^{s+1}e^{2t}(1-e^t)\Gamma(s+2)\sum_{k=0}^{\infty}e^{tk}\left[1-\frac{1}{(1+1/2^k)^{s+2}}\right]$$

and

$$\mathbf{Mellin}\left\{\frac{z^2}{2}e^{2t}e^{-z/2}; s\right\} = 2^{s+1}e^{2t}\Gamma(s+2).$$

The last piece gives

$$\mathbf{Mellin}\left\{\frac{z^2}{2}e^{2t}e^{-z}; s\right\} = \frac{1}{2}e^{2t}\int_0^{\infty}e^{-z}z^{s+1}\,dz$$
$$= \frac{1}{2}e^{2t}\Gamma(s+2).$$

Therefore, the second part on the right-hand side of (11) gives

$$\mathbf{Mellin}\left\{ze^{-z}e^{2t}\sum_{j=2}^{\infty}\frac{j\phi_{\delta_j}(t)}{j!}\left(\frac{z}{2}\right)^j; s\right\} = -2^{s+1}e^{2t}(1-e^t)\Gamma(s+2)$$
$$\times\sum_{k=0}^{\infty}e^{tk}\left[1-\frac{1}{(1+1/2^k)^{s+2}}\right]$$
$$+\,2^{s+1}e^{2t}\Gamma(s+2)-\frac{1}{2}e^{2t}\Gamma(s+2). \tag{15}$$

For the third and last part of the right-hand side of (11), we get that

$$\mathbf{Mellin}\left\{e^{2t}\frac{z^2}{4}(1-e^{-z}); s\right\} = \frac{1}{4}e^{2t}\int_0^{\infty}(1-e^{-z})z^{s+1}\,dz$$
$$= -\frac{1}{4}e^{2t}\Gamma(s+2). \tag{16}$$

After putting (14), (15) and (16) together and doing some algebra, equation (10) results in Lemma 1.

**Acknowledgment.** The authors are indebted to Dr. Ralph Neininger for sound advice on the scope of the contraction method in the presence of oscillation.

Department of Statistics
George Washington University
Washington, DC 20052
USA
e-mail: coschri@gwu.edu
e-mail: hosam@gwu.edu